\documentclass[12pt]{amsart}
 \usepackage{latexsym}
 \usepackage{amsmath}
 \usepackage{amsfonts}

\newtheorem{thm}{Theorem}
\newtheorem{corollary}[thm]{Corollary}
\newtheorem{example}[thm]{Example}
\newtheorem{lemma}[thm]{Lemma}
\newtheorem{proposition}[thm]{Proposition}

\newtheorem{conjecture}[thm]{Conjecture}
\theoremstyle{definition}
\newtheorem{definition}[thm]{Definition}
\theoremstyle{remark}

     \newcommand{\Z}{\mathbb{Z}}
     \newcommand{\N}{\mathbb{N}}
     \newcommand{\R}{\mathbb{R}}
     \renewcommand{\chi}{{\bf{1}}}
     \newcommand{\T}{\mathbb{T}}

     \newcommand{\calD}{\mathcal{D}}
     
     \newcommand{\cD}{\mathcal{D}}

     \newcommand{\La}{\Lambda}
     \newcommand{\la}{\lambda}

 \newcommand{\ceil}[1]{\lceil #1 \rceil}
 
 \newcommand{\floor}[1]{\lfloor #1 \rfloor}

\newcommand{\cal}{\mathcal}

\newcommand{\bt}{{\mathbb T}}

\newcommand{\bz}{{\mathbb Z}}

\newcommand{\lan}{\langle}
\newcommand{\ran}{\rangle}

\newcommand{\discr}{\operatorname{Discr}}

\newcommand{\calK}{{\cal K}}

\begin{document}

\title{Uniform partitions of frames of exponentials into Riesz sequences}

       \author{ Darrin Speegle}
\address      {  Department of Mathematics and Computer Science\\
        Saint Louis University\\
        221 North Grand Blvd., St. Louis, MO 63103, USA}
\email        { speegled@slu.edu}

\date{\today}
\footnotetext{\textit{Math Subject Classifications.} 42C15, 42C40}

\begin{abstract} 
The Feichtinger Conjecture, if true, would have as a corollary that for each set $E\subset \T$ and $\Lambda \subset \Z$, there is a partition $\Lambda_1,\ldots,\Lambda_N$ of $\Z$ such that for each $1\le i \le N$, $\{\exp(2\pi i x\lambda): \lambda \in \Lambda_i\}$ is a Riesz sequence.  In this paper,
sufficient conditions on sets $E\subset \T$ and $\Lambda\subset \R$ are given so that $\{\exp(2\pi i x\lambda) 1_E: \lambda \in \Lambda\}$ can be uniformly partitioned into Riesz sequences.     
\end{abstract}

\maketitle

\section{Introduction}   A frame is a collection of elements $\{e_i : i\in I\}$ in a
Hilbert space ${\mathcal H}$ such that there exist positive constants $A$ and $B$ such that for every $h\in \mathcal H,$
\[
A\|h\|^2 \le \sum_{i\in I} | \lan h, e_i \ran|^2 \le B\|h\|^2.
\]
A frame $\{e_i: i\in I\}$ is {\it bounded} if 
\[
\inf_{i\in I} \|e_i\| > 0.
\]
(Note that it is automatic that $\sup_{i\in I} \|e_i\| < \infty$.)  A sequence $\{e_i: i\in I\}$ is said to be a Riesz
sequence if it is a Riesz basis for its closed linear span, i.e.,
there exist $K_1, K_2>0$ such that for every finite family of scalars $\{a_i: i\in I\}$
\[
K_1 \sum_{i\in I} |a_i|^2 \le \bigg\|\sum_{i\in I} a_i e_i \bigg\|^2 \le K_2 \sum_{i\in I} |a_i|^2.
\]
Note that for Riesz basic sequences we always have
\[
0 < \inf_{i\in I}\|e_i\| \le \sup_{i\in I} \|e_i\| < \infty.
\]
With this notation, one can state the Feichtinger conjecture: 
\begin{conjecture}\label{feicon} (Feichtinger) Every bounded frame can be written as the finite disjoint union of Riesz basic sequences.
\end{conjecture}

A sort of converse to the Feichtinger is true; namely, the finite union of Riesz basic sequences is a frame for its closed linear span.

The Feichtinger conjecture has been shown to be related to several famous open problems in analysis, and as such, is receiving a fair amount of recent interest \cite{BS, CCLV, CT, Gr}.  Of particular interest is the special case of the Feichtinger conjecture for frames of exponentials, e.g. frames of the form $\{\exp(2\pi i \lambda x)1_E: \lambda \in \Lambda\}$.   The best positive result so far is in \cite{BT}, where it is shown that in the case $\Lambda = \Z$, frames of exponentials always contain a Riesz sequence with positive Beurling density.  

The most natural type of partition of a set indexed by the integers would be a uniform partition, so it is natural to ask which frames can be uniformly partitioned into Riesz sequences.  We state this formally as a definition.

\begin{definition}\label{defone}  Let $\Lambda = \{ \cdots < \lambda_{-1} < \lambda_0 < \lambda_1 < \cdots\} \subset \R$.  We say that $\{e_{\lambda_k}: k\in \Z\}$ can be uniformly partitioned into Riesz sequences if there exists an $N$ such that for $1\le J \le N$, $\{e_{\lambda_{mN + J}}: m\in \Z\}$ is a Riesz sequence.  
\end{definition}

Gr\"ochenig \cite{Gr} showed that if a frame is intrinsically localized, than it can be uniformly partitioned into Riesz sequences.  Bownik and the author \cite{BS} observed that if $E$ contains an interval a.e., then $\{\exp(2\pi i x n) 1_E: n\in \Z\}$ can be uniformly partitioned into Riesz sequences.  Halpern, Kaftal and Weiss \cite{HKW} showed that if $\phi \in L^\infty(\T)$ is Riemann integrable, then the associated Laurent operator $L_\phi$ can be uniformly paved.  Moreover, to date the primary negative evidence offered against the Feichtinger conjecture is two constructions of frames which can not be uniformly partitioned into Riesz sequences in a strong way, see \cite[Theorem 4.13]{BS} and \cite[Theorem 5.4 (b)]{HKW}.

In this paper, we provide a sufficient condition on the pair $(E, \Lambda)$, where $E\subset \T$ and $\Lambda \subset \R$, such that $\{\exp(2\pi i \lambda x)1_E: \lambda \in \Lambda\}$ can be uniformly partitioned into Riesz sequences.   Of particular interest is the application of this condition to the example considered in \cite{BT, HKW}, which shows that, perhaps, the proposed counterexample to the paving conjecture in \cite{HKW}, which was proven not to be a counterexample in \cite{BT}, was not optimally chosen.   See example \ref{example1} for details.

\section{Beurling dimension}
In this section, we recall some facts about the Beurling dimension of a subset of $\R^d$, though we will be concerned only with subsets of $\R$.    For $h > 0$, we let $Q$ denote the cube  $[-1,1]^d$ and let $Q_h$ be the dilation of $Q$ by the factor of $h$:
\[
Q_h = hQ = [-h, h]^d.
\]
For any $x = (x_1, \ldots, x_d) \in \R^d$ we let $Q_h(x)$ be the set $Q_h$ translated
in such a way so that it is ``centered'' at 
$x$, i.e.,
\[
Q_h(x)  = \prod_{i=1}^d [x_i - h, x_i + h].
\]
Employing these notions we will first define a generalization of Beurling density. 

\begin{definition}
Let $\La \subset \R^d$ and $r > 0$. Then the \emph{lower 
Beurling density of ${\La}$ with respect to $r$} is defined by
\[
{\mathcal{ D}}^-_r({{\La}}) = \liminf_{h\to\infty} \inf_{x \in \R^d}
\frac{\#({\La}\cap Q_h(x))}{h^r},
\]
and the \emph{upper Beurling density of ${\La}$ with 
respect to $r$} is defined by
\[
{\mathcal { D}}^+_r({\La}) = \limsup_{h\to\infty} \sup_{x \in \R^d}
\frac{\#({\La}\cap Q_h(x))}{h^r}.
\]
If $\cD^-_r(\La) = \cD^+_r(\La)$, then we say that $\La$ has \emph{uniform Beurling 
density with respect to $r$} and denote this density by $\cD_r(\La)$. 
\end{definition}

With this definition in hand, we can define the upper and lower Beurling dimensions of subsets of $\R^d$.
 
\begin{definition}
Let $\La \subset \R^d$. Then the \emph{lower dimension of ${\La}\subset \R^d$} is
defined by
\[
\dim^-({\La}) = \inf\, \{ r>0: {\mathcal { D}}_r^- ({\La}) < \infty \}
\]
and the \emph{upper dimension of ${\La}\subset \R^d$} is 
\[
\dim^+({\La}) = \sup\, \{ r>0: {\mathcal{ D}}_r^+ ({\La}) > 0 \}.
\]
When these two quantities are equal, we refer to the \emph{Beurling dimension} 
of ${\La}$, and we denote it by $\dim ({\La})$.
\end{definition}
 
 We note here that the upper Beurling dimension is a base point independent version of the {\it upper mass dimension} considered in \cite{BT89}, see \cite{CKS, CKS2} for details.  One result on Beurling dimension that we will use in the main section of this paper is the following \cite{CKS, CKS2}.
 
 \begin{proposition}
\label{p4}
Let ${\La} \subset \R^d$. 
\begin{enumerate}
\item The following conditions are equivalent.
\begin{enumerate}
\item $\calD^+_d({\La}) < \infty$.
\item There exists some $h > 0$ such that
$\sup_{x \in \R^d} \# ({\La} \cap Q_{h}(x)) < \infty$.
\item For all $h > 0$,
$\sup_{x \in \R^d} \# ({\La} \cap Q_{h}(x)) < \infty$.
\item ${\La}$ is relatively uniformly discrete.
\item For all $h > 0$, $\sup_{x \in \R^d} \#\{\la \in \La : x \in Q_h(\la)\} < \infty$.
\end{enumerate}
\item Also the following conditions are equivalent.
\begin{enumerate}
\item $\calD^-_d({\La}) > 0$.
\item There exists some $h > 0$ such that
$\inf_{x \in \R^d} \# ({\La} \cap Q_{h}(x)) > 0$.
\item $\La$ contains a subsequence of positive uniform density.
\item There exists some $h > 0$ such that ${\La}$ is $h$-dense.
\end{enumerate}
\end{enumerate}
\end{proposition}

\section{Results}

For $E$ a measurable subset of $\T$ and $\Lambda \subset \R$, we will say that $(E, \Lambda)$ can be uniformly partitioned into Riesz sequences if the frame  $\{\exp(2\pi i \lambda x) 1_E: \lambda \in \Lambda\}$ can be uniformly partitioned into Riesz sequences, see Definition \ref{defone}.  Characterizing sets $E$ such that $(E, \Z)$ can be uniformly partitioned into Riesz sequences is related (but not equivalent) to characterizing the Laurent operators which can be uniformly paved, which was characterized in \cite{HKW} as those Laurent operators whose symbol is Riemann integrable.  
In this article, we consider general sets $\Lambda \subset \R$.

Our main tool is a theorem due to Montgomery and Vaughan: 
\cite[Theorem 1, Chapter 7]{Mon},
\cite{MV74}.

\begin{thm}\label{mv} Suppose that $\lambda_1,\ldots, \lambda_N$ are distinct real numbers, and suppose that $\delta > 0$
is chosen so that $|\lambda_n - \lambda_m| \ge \delta$ whenever $n\not= m$.  Then, for any coefficients $a_1,\ldots, a_N$, and any
$T > 0$,
\begin{equation} \label{monriesz}
\bigl(T - 1/\delta\bigr) \sum_{n=1}^N |a_n|^2 \le \int_0^T \biggl| \sum_{n = 1}^N a_n e^{2\pi i \lambda_n t} \biggr|^2\,
dt \le \bigl(T + 1/\delta\bigr) \sum_{n=1}^N |a_n|^2.
\end{equation}
\end{thm}

We will also need two lemmas concerning partitioning subsets of $\R$ with finite upper dimension.

\begin{lemma}\label{l1}  Let $\Lambda = \{\cdots < \lambda_{-1} < \lambda_0 < \lambda_1 < \cdots\} \subset \R$ such that $\dim^+(\Lambda) \le 1$, and $K\in \N$.  There exists $N\in \N$ such that whenever $|i - j| > N$, $|\lambda_i - \lambda_j| > K$. 
\end{lemma}

\begin{proof}  This is just a restatement of Proposition~\ref{p4}.
\end{proof}

For $\Lambda \subset \R$, $\alpha \ge 0$ and $r> 0$, we define $D_{\alpha, \Lambda}^{+}(r) = \sup\{\frac {\#(\Lambda \cap Q_r(x)}{r^{\alpha}} : x\in \R\}$.

\begin{lemma}\label{l2}  Let $\Lambda = \{\cdots < \lambda_{-1} < \lambda_0 < \lambda_1 < \cdots\} \subset \R$ be such that $\dim^+(\Lambda) < \beta \le 1$,  and $\epsilon > 0$.  There exist $N, R\in \N$ such that for each $1\le j\le N$ and $r\ge R$, $D_{\beta, \Lambda_j(N)}^+(r) \le 2 R^{-\beta} + \epsilon$, where $\Lambda_j(N) = \{\lambda_{mN + j}: m\in \Z\}$.  
\end{lemma}

\begin{proof}  Choose $R$ such that for $r \ge R$, 
\[
\sup_{x\in \R} \frac{\#(\Lambda \cap Q_r(x))}{r^\beta} \le \sup_{x\in \R} \frac {\#(\Lambda \cap Q_R(x))} {R^\beta} + \epsilon < \infty.
\]
Choose $N = \sup_{x\in \R} \#(\Lambda \cap Q_R(x))$.  Now, fix $r\ge R$ and $1\le j \le N$.  Then, 
\begin{eqnarray*}
D^+_{\beta, \Lambda_j(N)}(r) &=& \sup_{x\in \R} \frac{\#(\Lambda_j(N) \cap Q_r(x))}{r^\beta}\\
&\le& \sup_{x\in \R} \frac {\#(\Lambda \cap Q_r(x))/N + 1}{r^\beta}\\
&\le& \sup_{x\in \R}\frac {\#(\Lambda \cap Q_R(x))}{R^\beta N} + \epsilon + R^{-\beta}\\
&\le& 2R^{-\beta} + \epsilon.
\end{eqnarray*}


\end{proof}

Using Theorem \ref{mv}, we obtain the following improvement to Lemma 5.1 in \cite{CCK}.

\begin{lemma}\label{beter}  Suppose $\Lambda \subset \R$.  If $I$ is an interval contained in $\T$, then for any sequence of numbers $\{a_\lambda\}_{\lambda \in \Lambda}$ in $\ell^2$, 
\begin{equation}
\int_I  |\sum a_\lambda e^{2\pi i \lambda \xi} |^2 \le 2 \ell(I) D^+_{0, \Lambda}(\ell(I)^{-1})
\end{equation}
\end{lemma}

\begin{proof} First note that the interval $[0,T]$ in Theorem \ref{mv} can be replaced by any interval of length $T$, which we set to be $\ell(I)$.   We can partition $\Lambda$ into $D^+_{0, \Lambda}(\ell(I)^{-1}) $ subsets $\Lambda_i$ such that if $\lambda_1, \lambda_2$ are in $\Lambda_{i}$, then $|\lambda_1 - \lambda_2| > \ell(I)^{-1}$.  It follows that 
\begin{eqnarray*}
\bigl(\int_I  |\sum a_\lambda e^{2\pi i \lambda \xi} |^2 \bigr)^{1/2} &\le& \sum_{i} \bigl(\int_I  |\sum_{\lambda \in \Lambda_{i}} a_\lambda e^{2\pi i \lambda \xi} |^2 \bigr)^{1/2} \\
&\le& \sum_{i} (2 \ell(I))^{1/2} \bigl(\sum_{\lambda \in \Lambda_{i}} |a_\lambda|^2\bigr)^{1/2} \\
&\le& (2 \ell(I))^{1/2} ( D^+_{0, \Lambda}(\ell(I)^{-1}) )^{1/2} \bigl(\sum_{i}\sum_{\lambda \in \Lambda_{i}} |a_\lambda|^2\bigr)^{1/2}\\
&=& (2\ell(I))^{1/2} (D^+_{0, \Lambda}(\ell(I)^{-1}))^{1/2} \bigl(\sum_{\lambda \in \Lambda} |a_\lambda|^2\bigr)^{1/2},
 \end{eqnarray*}
where the second inequality is from Theorem \ref{mv} and the third inequality is the generalized mean inequality with $D^+_{0, \Lambda}(\ell(I)^{-1})$ terms.
\end{proof}


We recall for motivation of the hypotheses in the following theorem the definition of the essential $\alpha$-Hausdorff measure of a set $E$ --- $H_{\alpha}(E) = \inf\{\sum \ell(I_{n})^\alpha: E \subset \cup_{n=1}^{\infty} I_{n} \cup J$, $|J| = 0\}$.  The following theorem is related to the essential $\alpha$-Hausdorff measure of $E$ in Corollary \ref{hausvers}.

\begin{thm} \label{maintheorem} Let $E\subset \T$ be measurable, $\Lambda = \{< \cdots < \lambda_1 < \lambda_0 < \lambda_1 < \cdots\}  \subset \R$ and $0 < \alpha < 1$.  If there exists a sequence of intervals $\{E_n:n\in \N\}$ of nonincreasing length, an integer $Z$ and $0 < \alpha < 1$ such that
\begin{enumerate}
\item  $\cup_{n=1}^\infty E_n \supset E$, and
\item  $\sum_{n = 1}^Z |E_n| + \sum_{n = Z + 1}^\infty |E_n|^{\alpha} < 1$,
\end{enumerate}
and $\dim^+(\Lambda) < 1 - \alpha$, then $(E, \Lambda)$ can be uniformly partitioned into Riesz sequences.
\end{thm}

\begin{proof}  Define $F =  \cup_{n=1}^\infty E_n$, and note that $|F| < |\T|$, and normalize $|\T| = 1$.  Choose $\epsilon > 0$ satisfying
\begin{enumerate}
\item $\frac 34 |\T| + \frac 14 |F| < 1 - \epsilon$ and
\item $|F| + 2\epsilon < \frac 12 {|\T| + |F|}.$
\end{enumerate}
Choose $M \ge Z \in \N$ such that 
\begin{enumerate}
\item $|E_M|^{1 - \alpha} < \frac {1}{4 \sum_{n = M + 1}^\infty |E_n|^\alpha}$, 
\item $4\bigl(|E_M|^{1 - \alpha} + \epsilon\bigr) \sum_{n = M + 1}^\infty |E_n|^\alpha < \epsilon$, and
\item $|E_M|^{-1} > R$, where $R$ is chosen from Lemma \ref{l2} with $\beta = 1- \alpha$ and $\epsilon$ as above.
\end{enumerate}

Choose $K \in \N$ such that $M/K < \epsilon$.   By Lemma \ref{l1}, there exists  $L \in \N$ such that $|\lambda_j - \lambda_k| > K$ whenever $|j - k| > L$.   By Lemma \ref{l2}, there exists  $J\in \N$ such that $D_{1 - \alpha, \Lambda_i(J)}^+(r) \le 2R^{\alpha - 1} + \epsilon \le  2|E_M|^{1 - \alpha} + \epsilon$ for all $r\ge |E_M|^{-1}$ and $1\le i \le J$.  Finally, let $N$ be the larger of $J$ and $L$.



Fix $1\le j \le N$, and $\{a_\lambda: \lambda \in \ell^2( \Lambda_j(N))\}$.  We compute
\begin{eqnarray*}
\sum_{n=1}^\infty \int_{E_n} |\sum_{\lambda\in \Lambda_j(N)} a_\lambda e^{2\pi i \lambda \xi}|^2 &=& 
\sum_{n = 1}^M\int_{E_n} |\sum_{\lambda\in \Lambda_j(N)} a_\lambda e^{2\pi i \lambda \xi}|^2  \\
&+&  \sum_{n = M + 1}^\infty    \int_{E_n} |\sum_{\lambda\in \Lambda_j(N)} a_\lambda e^{2\pi i \lambda \xi}|^2\ := S_1 + S_2.
\end{eqnarray*}

Moreover,  by Theorem \ref{mv} and our choice of $K$,
\begin{eqnarray}\label{S1}
S_1 \le  \biggl(\sum_{n = 1}^M \bigl(|E_n| + 1/K\bigr)\biggr)  \sum_{\lambda \in \Lambda_j(N)} |a_\lambda|^2 \le \bigl(|F| + \epsilon\bigr)  \sum_{\lambda \in \Lambda_j(N)} |a_\lambda|^2 .
\end{eqnarray}

We also have that

\begin{eqnarray*} 
S_2 &\le& 2\sum_{n=M + 1}^\infty \ell(E_n) D_{0, \Lambda_j(N)}^+((\ell(E_n)^{-1}) \sum_{\lambda\in \Lambda_j(N)} |a_\lambda|^2\\
&=& 2\sum_{n=M + 1}^\infty \ell(E_n) \ell(E_n)^{\alpha - 1} D_{\alpha, \Lambda_j(N)}^+((\ell(E_n)^{-1}) \sum_{\lambda\in \Lambda_j(N)} |a_\lambda|^2\\
&\le& 2\sum_{n=M + 1}^\infty \ell(E_n)^\alpha (2|E_M|^{1 - \alpha} + \epsilon) \sum_{\lambda\in \Lambda_j(N)} |a_\lambda|^2\\
&\le&4 (|E_M|^{1 - \alpha} + \epsilon) \sum_{n=M + 1}^\infty \ell(E_n)^\alpha \sum_{\lambda\in \Lambda_j(N)} |a_\lambda|^2 .
\end{eqnarray*}

In particular, 
\begin{equation}\label{S2}
S_2 < \epsilon \sum_{\lambda \in \Lambda_j(N)} |a_\lambda|^2.
\end{equation}

Therefore, combining (\ref{S1}) and (\ref{S2}), we get that 
\begin{eqnarray*}
\sum_{n=1}^\infty \int_{E_n} |\sum_{\lambda\in \Lambda_j(N)} a_\lambda e^{2\pi i \lambda \xi}|^2 &\le&
\bigl(|F| + 2\epsilon) \sum_{\lambda \in \Lambda_j(N)} |a_\lambda|^2\\
& <& (|\T| + |F|)/2 \sum_{\lambda \in \Lambda_j(N)}|a_\lambda|^2.
\end{eqnarray*}
Therefore, since 
\begin{eqnarray*}
\frac{3|\T| + |F|}{4} \sum_{\lambda \in \Lambda_j(N)} |a_\lambda|^2 &\le& (1 - 1/K)\sum_{\lambda \in \Lambda_j(N)} |a_\lambda|^2 \\ 
&\le& \int_{\T} \bigl|\sum_{\lambda\in \Lambda_j(N)} a_\lambda e^{2\pi i \lambda \xi}\bigr|^2 \, d\xi\\
&=& \bigl( \int_{E^c} + \int_{E} \bigr) \bigl|\sum_{\lambda\in \Lambda_j(N)} a_\lambda e^{2\pi i \lambda \xi}\bigr|^2 \, d\xi \\
&\le& \bigl(\int_{E^c} + \int_{\cup E_n}\bigr)  \bigl|\sum_{\lambda\in \Lambda_j(N)} a_\lambda e^{2\pi i \lambda \xi}\bigr|^2\, d\xi\\
&\le&  \int_{E^c}  \bigl|\sum_{\lambda\in \Lambda_j(N)} a_\lambda e^{2\pi i \lambda \xi}\bigr|^2\,d\xi + (|\T| + |F|)/2 \sum_{\lambda\in \Lambda_j(N)} |a_\lambda|^2,
\end{eqnarray*}
it follows that 
\[
(|\T| - |F|)/4 \sum_{\lambda \in \Lambda_j(N)} |a_\lambda|^2 <  \int_{E^c}  |\sum_{\lambda\in \Lambda_j(N)} a_\lambda e^{2\pi i \lambda \xi}|^2,
\]
i.e.~$\{1_{E^c} e^{2\pi i \lambda \xi}: \lambda \in \Lambda\}$ has  a lower Riesz bound.  Since each $\Lambda_j(N)$ is separated, it also has an upper Riesz bound, which completes the proof.  
\end{proof}

\begin{corollary}\label{hausvers}   Let $E\subset \T$ be measurable, $\Lambda = \{< \cdots < \lambda_1 < \lambda_0 < \lambda_1 < \cdots\}  \subset \R$ and $0 < \alpha < 1$.  If $H_\alpha(E) < 1$ and $\dim^+(\Lambda) < 1 - \alpha$, then $(E, \Lambda)$ can be uniformly partitioned into Riesz sequences.
\end{corollary}





\begin{example} \label{example1}  An example considered by Bourgain-Tzafriri \cite{BT} and Halpern-Kaftal-Weiss \cite{HKW}. 
{\rm Let $\{r_n:n\in \N\}$ be a partition of the rational numbers in $\T$.   For each $n$, let $E_n$ be an interval centered at $r_n$ with length $|E_n| < 2^{-n}|\T|$.   Let $\phi$ be the indicator function supported on $E = \cup_{n = 1}^\infty E_n$.  In \cite{HKW}, it was shown that the Laurent operator with symbol $\phi$,  $L_\phi$, cannot be uniformly paved (in fact, something stronger was shown), while in \cite{BT}, it was shown that $L_\phi$ can still be paved.  In particular, this implies that there is a partition of the integers $\Lambda_1,\ldots,\Lambda_N$ such that  for each $1\le i \le N$, $\{1_{E^c} \exp(\lambda): \lambda \in \Lambda_i\}$ is a Riesz sequence.  (Note that $(E^c, \Z)$ cannot be uniformly partitioned into Riesz sequences.)  However, by Theorem \ref{maintheorem}, whenever $\dim^+(\Lambda) < 1$, $(E, \Lambda)$ can be uniformly partitioned into Riesz sequences.  So, perhaps a better candidate for a counterexample to the Feichtinger Conjecture would have been to have a series $\sum a_n$ that converges more slowly to a number less than 1 and to let $|E_n| < a_n |\T|$.}
\end{example}

We end this paper by presenting an example of  a set $E\subset \T$ and a set $\Lambda\subset \Z$ such that $\dim^+(\Lambda) < 1$ yet $(E, \Lambda)$ can not be uniformly partitioned into Riesz sequences.  We begin by recalling the following theorem.

\begin{thm}\label{strangeset} \cite{BS} There exists a set $E \subset \T$ such 
that whenever $\calK \subset \bz$ is such that for all $\delta > 0$, there exist $M,N \ell\in \bz$ such that 
\begin{enumerate}
\item\label{een} $\ell N^{-1/2} \log^3 N < \delta$, and
\item\label{twee} $\{M, M + \ell, \ldots, M + N\ell\} \subset \calK$,
\end{enumerate}
then $\exp(\lambda): \lambda \in \calK\}$ is not a Riesz sequence in $L^2(E)$.
\end{thm}

\begin{corollary}\label{Easycor}  There exists a set $E \subset \T$ such that for each $1 > \beta > 2/3$ there is a set $\Lambda \subset \Z$ such that $\dim^+(\Lambda) = \beta$ and $(E, \Lambda)$ cannot be uniformly partitioned into Riesz sequences.
\end{corollary}

\begin{proof}  Let $E$ be the set guaranteed to exist from Theorem \ref{strangeset}.  Let $2/3 < \beta < 1$.  Let $\gamma = \frac{1-\beta}\beta$.  Define $\Lambda = \cup_{j = 1}^\infty \{Q_j + \alpha \ceil{j^\gamma}: 0\le \alpha \le j\}$, where the $Q_j$'s are chosen to be some rapidly increasing sequence such as $2^{2^j}$.     It follows that $\sup_{x\in \R} |\# (\Lambda \cap Q_{j^{\gamma + 1}}(x))| \approx j$, and so $\dim^+(\Lambda) = \frac 1{1+\gamma} = \beta$.  

Now, let $N$ be a positive integer and write $\Lambda = \{\lambda_n:n\in \N\}$ where $\lambda_i < \lambda_j$ when $i < j$.  Let $\Lambda_N = \{\lambda_{Nn}:n\in \N\}$.  We show that $\Lambda_N$ satisfies (\ref{een}) and (\ref{twee}) of the hypotheses of Theorem \ref{strangeset}, which completes the proof of this corollary.  Let $\delta > 0$.  For each $j$, let $k_j$ be the smallest nonnegative inteer such that $Q_j + k_j \ceil{j^\gamma} \in \Lambda_N$.  We have that 
\[
\{Q_j + (\alpha N + k_j)\ceil{j^\gamma} : 0\le \alpha < \floor{j/N}\}\subset \Lambda_N.
\]
Since $\gamma < 1/2$, we can find $j$ such that 
\[
\floor{j/N}^{-1/2} N\ceil{j^\gamma} \log^3 \floor{j/N}  < \delta,
\]
which finishes the proof.
\end{proof}


\section*{Acknowledgments}

The author wishes to thank Dick Gundy, whose insightful questions at a seminar at Washington University led to some improvements to this paper, and Wojtek Czaja and Gitta Kutyniok for their encouragement to pursue this line of thought.  The author is partially supported by NSF-DMS 0354957.


\begin{thebibliography}{99}





 
\bibitem{BT2}ÊJ.~Bourgain and L.~Tzafriri, Invertibility of ``large" submatrices with applications to the geometry of Banach
spaces and harmonic analysis. {\sl Israel J. Math.} {\bf 57} (1987), no. 2, 137--224.

\bibitem{BT} J.~Bourgain and L.~Tzafriri, On a problem of Kadison and Singer. {\sl J. Reine Angew. Math.} {\bf 420} (1991),
1--43.


\bibitem{BS} M. Bownik and D. Speegle, The Feichtinger conjecture for wavelet frames, Gabor frames, and frames of translates, {\sl Can. J. Math.} {\bf 58} (2006), no. 6, 1121--1143.

  \bibitem{BT89}
  \rm{M.~T.~Barlow and S.~J.~Taylor,}
  \sl{Fractional dimension of sets in discrete spaces,}
  \rm{J. Phys. A: Math. Gen. {\bf 22} (1989), 2621--2626.}

\bibitem{CCK}
P. Casazza, O. Christensen, N. Kalton,  
Frames of translates.  
{\sl Collect. Math.} {\bf 52} (2001), no.~1, 35--54.

\bibitem{CCLV}
P. Casazza, O. Christensen, A. Lindner and R. Vershynin,
Frames and the Feichtinger conjecture,  {\sl Proc. Amer. Math. Soc.} {\bf 133}  (2005),  no. 4, 1025--1033. 

\bibitem {CT}
P, Casazza and J. Tremain, The Kadison-Singer problem in mathematics and engineering,  {\sl Proc. Natl. Acad. Sci. USA}  {\bf 103}  (2006),  no. 7, 2032--2039.

 

  \bibitem{CKS}
  \rm{W.~Czaja, G.~Kutyniok, D.~Speegle,} The geometry of sets of parameters of wave packet frames, 
{\sl Appl. Comput. Harmon. Anal.} {\bf 20}  (2006),  no. 1, 108--125.
  

\bibitem{CKS2}
  \rm{W.~Czaja, G.~Kutyniok, D.~Speegle,} 
  {\sl  Beurling dimension of Gabor pseudoframes of affine subspaces,}
 J. Fourier Anal. Appl. to appear. 


 \bibitem{Gr}
K. Gr{\"o}chenig,
Localized frames are finite unions of {R}iesz sequences, {\sl Adv. 
Comput. Math.} {\bf 18}
(2003), 149--157.

\bibitem{HKW} H.~Halpern, V.~Kaftal and G.~Weiss, Matrix pavings and Laurent operators. {\sl J.~Operator Theory} {\bf 16}
(1986), no.~2, 355--374.


  

\bibitem{Mon}
H. L. Montgomery, Ten lectures on the interface between analytic number theory and harmonic analysis. CBMS Regional Conference Series in Mathematics, 84. {\sl American Mathematical Society, Providence, RI,} 1994. xiv+220 pp.

\bibitem{MV74}
H.~L.~Montgomery and R.~C.~Vaughan,
Hilbert's Inequality.
{\sl J. London Math. Soc.} (2) {\bf 8} (1974), 73--82.



\end{thebibliography}
\end{document}